**Francesco Trimarchi**

**Rational points on elliptic curves and representations of rational numbers as the product of two rational factors**

## §1. *Introduction and summary*

In this paper we prove two propositions concerning:

i) the representations of rational numbers as the product of two rational factors;

ii) the related properties of elliptic curves such that the cubic has rational roots.

The first proposition (Iso-additive Representations Theorem, shortly IRT) states that any pair $(m, n)$ of non-zero and distinct rational numbers may have, at most, four representations $m = m_1 m_2$ and $n = n_1 n_2$ such that $m_1 + m_2 = n_1 + n_2$ where $m_1, m_2, n_1, n_2 \in \mathbb{Q}_{\neq 0}$ and $m_1 m_2 \neq n_1 n_2$ (Iso-additive Representations, shortly IR).

Here we consider only the representations essentially different, i.e. we count only once the two representations of the same number that differs only for the order in which its factors are written; instead, we keep distinct the two representations of the same number such that one is obtained on the other by changing the sign of both factors[1].

Although these representations do not seem to ever been studied, they are implicitly involved in many problems of number theory. For example,

i) if $m$ and $n$ are natural strictly positive numbers such that $n = m + 2$, then the pair $(m, -n)$ has the IR $m = 1 \times m$ and $-n = -1 \times n$ such that the sum of factors is $m + 1$; therefore, the twin primes conjecture can be restated as follows: if $m$ and $n$ are natural distinct primes then there are infinitely many pairs of the type $(m, -n)$ having an IR;

ii) let $(a, b, c)$ a Pythagorean triple such that $a^2 + b^2 = c^2$ and $a = 2pq$, $b = p^2 - q^2$, $c = p^2 + q^2$ its primitive representation; then

- each pair $(a^2, b^2)$, $(a^2, c^2)$ and $(b, c^2)$ has the IR $a^2 = (2p^2) \times (2q^2)$, $b^2 = (p + q)^2 \times (p - q)^2$ and $c^2 = (p^2 + q^2) \times (p^2 + q^2)$ such that the sums of factors is always $2(p^2 + q^2)$,

---

[1] For example, let $m = 30$ and $n = 42$. The pair $(30, 42)$ has the following IR:

- $30 = 3 \times 10 = (-3) \times (-10)$ and $42 = 6 \times 7 = (-6) \times (-7)$,

- $30 = 2 \times 15 = (-2) \times (-15)$ and $42 = 3 \times 14 = (-3) \times (-14)$

such that the sum of factors is respectively $\pm 13$ and $\pm 17$. The IRT ensures that there are no other representations of this type. Unlike our example, the proof of the existence of two (or four) IR could in general be a prohibitive undertaking by means of a direct calculation. Indeed, any integer (or rational) number has infinitely many representations as a product of two rational factors or as a product of an integer and a rational factor. To find a pair of integers having four IR as a product of two integer is a lucky case.

- the pairs $(a^2, -b^2)$ and $(-a^2, b^2)$ have an IR: $a^2 = (2pq) \times (2pq)$ and $-b^2 = (p+q)^2 \times (-(p-q)^2)$ in the former case; $-a^2 = (2p^2) \times (-2q^2)$ and $b^2 = (p^2 - q^2) \times (p^2 - q^2)$ in the latter; here the sums of factors are, respectively, $4pq$ and $2(p^2 + q^2)^2$;

iii) let $p = 2m + 1$ and $q = 2n + 1$ distinct odd primes; as a consequence of Fermat Theorem there are two strictly positive integers $k$ and $h$ such that $1 - x^{2m} = -kp$ and $1 - x^{2n} = -hq^3$; then, for any pair of odd primes $(p, q)$ there is a pair $(-kp, -hq)$ having an IR such that the sum of factors is two: $-kp = (1 + x^m) \times (1 - x^m)$ and $-hq = (1 + x^n) \times (1 - x^n)$.

In Ch. 2 we give a constructive proof which requires only elementary techniques. Namely, we prove that that for every pair $(m, n)$ of non-zero and distinct rational numbers the equation

$$(*) \quad X + \frac{m}{X} = Y + \frac{n}{Y}$$

has at most a finite number of rational solutions[4]; we also give the parametric formulas of any pair $(m, n)$ having at least one IR, as well the formulas of their factors.

Solving this equation in rational (or integer) numbers is a typical Diophantine problem, that can be stated in rhetorical form in several equivalent ways; for example (in brackets the corresponding Diophantine equation to be solved):

i) *"find two numbers such that their ratio is equal to the ratio between their product diminished by a given number, and their product diminished by another given number"* $\left(\frac{X}{Y} = \frac{XY - m}{XY - n}\right)$;

ii) *"find two numbers such that the ratio between the first and the second is equal to the ratio between the square of the first increased by a given number, and the square of the second increased by another given number"* $\left(\frac{X}{Y} = \frac{X^2 + m}{Y^2 + n}\right)$.

---

[2] Since $-a^2 = -a \times a$, $-b^2 = -b \times b$ and $-c^2 = -c \times c$, also the pairs $(-a^2, -b^2)$, $(-a^2, -c^2)$ and $(-b^2, -c^2)$ have an IR such that the sum of factors is zero; however, this property holds for any triple of integers (or rational) numbers and not only for the Pytagorean ones.

[3] Davenport, pp. 35-37.

[4] In the following the Latin capital letters indicate variables while the lowercase letters indicate the parameters.

Here the "*two numbers*" to be "*find*" are the unknowns $X$ and $Y$, the "*given numbers*" are $m$ and $n$ and both equation above are obtained manipulating the equation (\*)[5].

No one should be surprised of finding these problems ever in the Arithmetica of Diophantus and/or in the works of Fermat and Bachet; however, in these works there seems to be no reference to them[6]; on the other hand, nobody seems to have studied the equation $X + \frac{m}{X} = Y + \frac{n}{Y}$.

The second proposition (Iso-additive Rational Points Theorem, shortly IRPT) immediately follows from the first and concerns the rational points of the elliptic curve $E$: $W^2 = U(U-m)(U-n)$, where $m \neq n \in \mathbb{Q}_{\neq 0}$[7].

In Ch. 3 we show (Iso-additive Rational Points Theorem, shortly IRPT) that if the equation (\*) has a rational solution $p = (x, y)$ for the variables $X$ and $Y$ where $x \neq y \in \mathbb{Q}_{\neq 0}$ - i.e. the pair $(m,n)$ has at least the IR $m = x \cdot \frac{m}{x}$ and $y = y \cdot \frac{n}{y}$ such that $x + \frac{m}{x} = y + \frac{n}{y}$, or equivalently $\frac{m}{x} - y = \frac{n}{y} - x$ - then $\left[\left(\frac{m}{x} - y\right)xy\right]^2 = xy(xy-m)(xy-n)$.

Therefore, by the transformation $XY \rightarrow U$ and $\left(\frac{m}{x} - Y\right)XY \rightarrow W$, or equivalently $\left(\frac{n}{Y} - X\right)XY \rightarrow W$, to any solution $p = (x, y)$ of the equation (\*) corresponds on $E$ a rational point $P = (u,w)$, where $u = xy$ and $w = \left(\frac{m}{x} - y\right)xy$.

This transformation defines a map from the set, say $\Gamma(m,n)$, whose elements are the rational solution of the equation (\*) and a set, say $\Omega(m,n)$ which is a sub-set of the group of rational points on $E$ and, as a consequence of IRT, may have up to sixteen distinct elements; under some conditions these points are of finite order and can be computed using the explicit formulas given in Ch. 3, instead of ordinary methods based on Nagell-Lutz Theorem.

The family of elliptic curves such that the non-zero roots of the cubic have at least an IR includes some interesting members, for example:

i) the curve $W^2 = U(U-1)(U+3)$. It has only eight rational point which are of finite order; it implies that there is no arithmetic sequence $\{a_1, a_2, a_3, a_4\}$ such that all terms are distinct rational squares. Note that the roots of the cubic have the IR $1 = 1 \times 1$ and $-3 = -1 \times 3$ such that the sum of factors is two; as we will show, IRT and IRPT give the explicit formulas for computing all the torsion points without computing the divisors of discriminant and check if they are the $X$-coordinates of some rational point on the elliptic curve;

---

[5] These problems are equivalent to a third that can be stated as follows (in brackets the corresponding system of two Diophantine equation to be solved): "*find three squares such that the first diminished by the second is a given number and the first diminished by the third is another given number*" ($T^2 - R^2 = m$ and $T^2 - S^2 = n$) [5]. Here the "*three numbers*" to be "*find*" are the unknowns $R$, $S$ and $T$ and the "*given numbers*" are still $m$ e $n$ (see note 8).

[6] De Fermat.

[7] Here we use the letters $U$ and $W$, rather than the standard ones ($X$ and $Y$) for distinguish this curve from the equation (\*). The assumption $m \neq n \in \mathbb{Q}_{\neq 0}$ ensures that: i) the curve is non-singular; ii) also the coefficients of the cubic are rational.

ii) the curve $W^2 = U(U - m^3)(U - n^3)$ such that $m^3 = m_1^3 m_2^3$, $n^3 = n_1^3 n_2^3$ and $m_1^3 + m_2^3 = n_1^3 + n_2^3$. It is well-known that there are infinitely many rational 4-tuple $(m_1^3, m_2^3, n_1^3, n_2^3)$ such that $m_1^3 + m_2^3 = n_1^3 + n_2^3$, representable by the formulas of Euler-Binet[8]; then, there are infinitely many pairs of cubes of the type $m^3 = m_1^3 m_2^3$ and $n^3 = n_1^3 n_2^3$ having at least one IR; for example $90^3 = 9^3 \cdot 10^3$ and $12^3 = 1^3 \cdot 12^3$. The sum of cubic factors of these numbers coincides with a famous number (*taxicab* 1729);

iii) the special case of the Frey curve $W^2 = U(U - m^{2k})(U + n^{2k})$. If the identity $m^{2k} + n^{2k} = l^{2k}$ were true for some integer $m, n, k > 0$[9], then the two non-zero roots of the cubic would have the IR $m^{2k} = m^k \cdot m^k$ and $-n^{2k} = (m^k + l^k) \cdot (m^k - l^k)$ such that the sum of factors is $2m^k$. We discuss some implications if this fact in Ch. 3.

## §2 *The Iso − additive Representations Theorem*

For every $m \neq n \in \mathbb{Q}_{\neq 0}$ let $\rho(m, n)$ the number of IR essentially different (as defined in note 2). Since $m_1 + m_2 = n_1 + n_2 \implies (-m_1) + (-m_2) = (-n_1) + (-n_2)$, $\rho(m, n)$ is zero or even and the IRT formally asserts:

*Theorem* 1 (*Iso − additive Representation Theorem − IRT*) $\forall\, m \neq n \in \mathbb{Q}_{\neq 0}$, $\rho(m, n) \in \{0,\ 2,\ 4\}$.

Proof. We deduce the IRT from the properties of the set, say $\Gamma$, of rational (non-trivial) solutions of the equation

[1] $X + \frac{M}{X} = Y + \frac{N}{Y}$

Clearly this equation has infinitely many rational solutions; indeed, if we set $X = x$ and $Y = y$ where $x, y \in \mathbb{Q}_{\neq 0}$, then the set of rational solutions for the variables $M$ and $N$ is the rational line $N = y(x - y) + \frac{y}{x} M$ so that, for arbitrary $M = m$, there are infinite 4-tuple $p = (x, y, m, n)$, where $n = y(x - y) + \frac{y}{x} m$, satisfying the equation [1] for the variables $X, Y, M, N$. Thus, there are infinitely many pairs $(m, n)$ having an IR.

---

On the other hand, if we set $M = N = a$ and $X = Y = b$, where $b \neq 0$, the equation [1] has infinitely many trivial solution $p = (b, b, a, a)$ for the variables $X, Y, M, N$, reflecting the fact that two equal numbers have the same (infinitely many) representations.

But our goal is to find what happens for fixed $M = m$ and $N = n$, when $m \neq n \in \mathbb{Q}_{\neq 0}$ and to prove that, for every fixed pair $(m, n)$ of non-zero distinct rational numbers, the number of essentially different solutions is at most four.

Then let $\Gamma(m, n) \subset \Gamma$ the sub-set of solutions of [1] such that $M = m$ and $N = n$, corresponding to the set of solutions of the equation

[2] $X + \frac{m}{X} = Y + \frac{n}{Y}$.

Consequently $\Gamma = \bigcup_{m \neq n \in \mathbb{Q}_{\neq 0}} \Gamma(m, n)$ and for given $m$ and $n$ two cases may occur:

i)    $\Gamma(m, n) = \emptyset$ ; the equation $X + \frac{m}{X} = Y + \frac{n}{Y}$ has not rational solutions and $m$ and $n$ have no IR;

ii)   $\Gamma(m, n) \neq \emptyset$ ; the equation $X + \frac{m}{X} - \left(Y + \frac{N}{Y}\right) = 0$ has at least one rational solution; then $\Gamma(m, n)$ contains at least a pair $p_1 = (x, y)$, $x \neq y \in \mathbb{Q}_{\neq 0}$, whose elements satisfy [2] respectively for the variables $X$ and $Y$ so that $m$ and $n$ have the IR $m = x \cdot \frac{m}{x}$ and $n = y \cdot \frac{n}{y}$ such that $x + \frac{m}{x} = y + \frac{n}{y}$.

In former case, $\Gamma(m, n)$ must contain in addition to $p_1$ any other pair that can be obtained from $p_1$ by substituting

i)    $x$ and/or $y$ with their cofactors (respectively $\frac{m}{x}$ and $\frac{n}{y}$),

ii)   both the factors of $m$ e $n$ that appear in $p_1$ - as well as in the pairs derived from $p_1$ through the substitutions defined in point i) - with their opposites[10].

---

[10] Note that $\forall\, m \neq n \in \mathbb{Q}_{\neq 0}\ \rho(m,n) > 0$ if, and only if, $m$ and $n$ are of the type $m = t^2 - r^2$ and $n = t^2 - s^2$, where $t \in \mathbb{Q}$ and $r \neq s \in \mathbb{Q}$. The first part of this statement is very easy to prove since $m$ and $n$ have the obvious IR $m = (t+r)(t-r)$ and $n = (t+s)(t-s)$ such that the sum of factors is $2t$. On the other hand, if $p = (m,n,x,y)$ is a non-trivial rational point on $X + \frac{M}{X} = Y + \frac{N}{Y}$ then the identities $m = x(z-x)$ and $n = y(z-y)$ hold, where $z = x + \frac{m}{x} (= y + \frac{n}{y})$ and $x \neq y$. By applying the formula $uw = \left(\frac{u+w}{2}\right)^2 - \left(\frac{u-w}{2}\right)^2$, we can represent $m$ and $n$ as a function of $x$, $y$ and $z$ in the following way: $m = \left(\frac{z}{2}\right)^2 - \left(\frac{z-2x}{2}\right)^2$, $n = \left(\frac{z}{2}\right)^2 - \left(\frac{z-2y}{2}\right)^2$. Thus, if we set $t = \frac{z}{2}$, $r = \frac{z-2x}{2}$ and $s = \frac{z-2y}{2}$, also the second part of the statement is proved∎. This is the basis for solving

Let $p(m,n) = \{\pm p_j, \; j = 1,2,3,4\} \subseteq \Gamma(m,n)$ the sub-set of these pairs[11]:

$p_1 = (x, \; y), \quad -p_1 = (-x, -y),$

$p_2 = (x, \; \frac{n}{y}), \quad -p_2 = \left(-x, -\frac{n}{y}\right),$

$p_3 = (\frac{m}{x}, \; y), \quad -p_3 = (-\frac{m}{x}, -y),$

$p_4 = (\frac{m}{x}, \; \frac{n}{y}), \quad -p_4 = (-\frac{m}{x}, \; -\frac{n}{y})$

If there was a further point $q_1 = (x', y') \in \Gamma(m,n)$ such that $x' \neq y'$, $x \neq x'$ and $y \neq y'$ then, as before, $\Gamma(m,n)$ should contain, in addition to $q_1$, any other pair that can be obtained from $q_1$ through the substitutions defined above, i.e. all the pairs of the sub-set $q(m,n) = \{\pm q_j, \; j = 1,2,3,4\}$ where

$q_1 = (x', \; y'), \quad -q_1 = (-x', -y'),$

$q_2 = (x', \; \frac{n}{y'}), \quad -q_2 = \left(-x', -\frac{n}{y''}\right),$

$q_3 = (\frac{m}{x'}, \; y'), \quad -q_3 = (-\frac{m}{x'}, -y'),$

$q_4 = (\frac{m}{x'}, \; \frac{n}{y'}), \quad -q_4 = (-\frac{m}{x'}, \; -\frac{n}{y'})$

In principle, there could be other pairs, say $r_1 = (x'', y'')$, $s_1 = (x''', y''')$, …, satisfying the equation [1] which in turn generate respectively the sub-sets $r(m,n)$, $s(m,n)$, …. However, only one pair of $p(m,n)$, $q(m,n)$, $r(m,n)$, $s(m,n)$, … is relevant for computing $\rho(m,n)$; the others are not essentially different because obtained from these pairs simply by switching any factor with its cofactor[12].

Thus, we can arbitrarily chose the pair $(p_1, -p_1)$ from $p(m,n)$, the pair $(q_1, -q_1)$ from $q(m,n)$, … to form the sub-set $\Gamma_0(m,n) \subseteq \Gamma(m,n)$ whose elements are the solutions essentially different of the equation [2].

---

the third Diophantine problem stated in Ch. 1 (see note 4) and for an alternative proof of IRT.

[11] In general $\#p = 8$; sometimes, however, the elements of $p$ are two-by-two coinciding; this occurs when $x$ and/or $y$ coincide with the respective cofactors, that is when $x = \frac{m}{x}$ ($\Longrightarrow p_1 = p_3$, $p_2 = p_4$ and $m = x^2$) and/or $y = \frac{n}{y}$ ($\Longrightarrow p_1 = p_2$, $p_3 = p_4$ and $n = y^2$); in this case $\#p = 4$; the simplest numerical example is the pair $(3, \; 4)$; a more general case is given by the pairs of the type $(-m^2, \; -n^2)$.

[12] However, in Ch 3. we will show all elements of $p(m,n)$ and eventually $q(m,n)$ - not only the pairs $(p_1, -p_1)$ and eventually $(q_1, -q_1)$ - are relevant because there is a map $p(m,n) \rightarrow \mathcal{P}(m,n)$ and eventually $q(m,n) \rightarrow \mathcal{Q}(m,n)$, where $\mathcal{P}(m,n)$ and $\mathcal{Q}(m,n)$ are in general distinct sub-sets of $E(\mathbb{Q})$ the set of rational point on the curve $W^2 = U(U-m)(U-n)$, where $m \neq n \in \mathbb{Q}_{\neq 0}$.

Since from the definition $\rho(m, n) \equiv \#\Gamma_0(m,n)$ the IRT can be restated as follows:

if $\Gamma_0(m,n) \neq \emptyset$ (i.e. $\#\Gamma_0(m,n) \neq 0$) then

    i)    either $\Gamma_0(m,n) = \{p_1, -p_1\}$ so that $\#\Gamma_0(m,n) = 2$

    ii)   or, at most, $\Gamma_0(m,n) = \{p_1, -p_1, q_1, -q_1\}$ so that $\#\Gamma_0(m,n) = 4$.

Now rewrite the equation [1] as follows

[3] $XY^2 - (X^2 + M)Y + NX = 0$

By applying elementary techniques we will obtain the parametric formulas for all non-trivial solutions of the equation [3] (i.e. the elements of $\Gamma$); this allows us to characterize any sub-set $\Gamma(m,n) \subset \Gamma$ and $\Gamma_0(m,n) \subset \Gamma(m,n)$ for every pair $(m,n)$ where $m \neq n \in \mathbb{Q}_{\neq 0}$.

Now the equation [3] can be treated as a polynomial equation of degree two in $Y$, with rational coefficient $X$, $-(X^2 + M)$, $NX$, roots $\frac{1}{2X}(X^2 + M \pm \sqrt{\Delta_Y})$ and discriminant $D_Y = X^4 + 2(M - 2N)X^2 + M^2$.

Therefore, the rational solutions of [3] coincides with rational points on the algebraic curve

[4] $\qquad Y = \frac{1}{2X}(X^2 + M \pm \sqrt{D_Y})$

and the necessary and sufficient condition for the existence of points of this type is $D_Y = T^2$, or equivalently $D_Y - T^2 = 0$, so that the [4] becomes

[5] $\qquad Y = \frac{1}{2X}(X^2 + M \pm T)$

Now $D_Y - T^2 = 0$ is equivalent to

[6] $\qquad X^4 + 2(M - 2N)X^2 + M^2 - T^2 = 0$

which in turn can be treated as a polynomial equation of degree four in $X$, with rational coefficients $1$, $0$, $2(M - 2N)$, $0$, $M^2 - T^2$ and roots $\pm\sqrt{2N - M \pm \sqrt{D_X}}$, where $D_X = 4N(N - M) + T^2$.

Therefore, the rational solutions of [6] coincides with rational points on the algebraic curve

[7] $\qquad X = \pm\sqrt{2N - M \pm \sqrt{D_X}}$

and the necessary and sufficient condition for the existence of points of this type is that at least one of the two equation

[8.1] $\qquad L_1^2 = 2N - M + \sqrt{D_X}$

[8.2]    $L_2^2 = 2N - M - \sqrt{D_X}$

is satisfied in rational numbers. But in turn [8.1]-[8.2] have rational solutions if, and only if, the equation $D_X = S^2$ has solutions of this type so that the [8.1]-[8.2] become

[9.1]    $L_1^2 = 2N - M + S$

[9.2]    $L_2^2 = 2N - M - S$

and, from [7], we have respectively

[10.1]    $X = \pm L_1,$

[10.2]    $X = \pm L_2$

Now, we rewrite the equation $D_X = S^2$ as follows

$$\frac{S+T}{2N} \cdot \frac{S-T}{2(N-M)} = 1$$

and obtain its (infinite) parametric solutions from the system

$$\begin{cases} \dfrac{S+T}{2N} = \dfrac{a}{b} \\ \dfrac{S-T}{2(N-M)} = \dfrac{b}{a} \end{cases}$$

where $a, b \in \mathbb{Z}_{\neq 0}$ are coprime.

Now rewrite $S$ and $T$ as a function of $N$ and $M$

[11.1]    $S = \frac{a^2+b^2}{ab} N - \frac{b}{a} M$

[11.2]    $T = \frac{a^2-b^2}{ab} N + \frac{b}{a} M$

and consider the following three distinct cases:

**Case** 1)  $a = b \in \mathbb{Z}_{\neq 0} \Rightarrow \frac{a}{b} = \frac{b}{a} = 1$

**Case** 2)  $a = -b \in \mathbb{Z}_{\neq 0} \Rightarrow \frac{a}{b} = \frac{b}{a} = -1$

**Case** 3)  $a \neq \pm b \in \mathbb{Z}_{\neq 0} \Rightarrow \frac{a}{b} \neq \pm \frac{b}{a}, \ \frac{a}{b} \neq \pm 1, \ \frac{b}{a} \neq \pm 1.$

In the following, we show that in any case, at most, $\rho(m, n) = 4$. Which demonstrates the IRT.

| **Case 1 (a = b)** |
|---|

The [11.1] and [11.2] become respectively

[12.1]    $S = 2N - M$

[12.2]    $T = M$

so that [9.1] and [9.2] are reduced, respectively, to

[13.1]    $L_1^2 = 2(2N - M) \Rightarrow X = \pm L_1$

[13.2]    $L_2^2 = 0 \Rightarrow X = 0$

Excluding the singular case $X = 0$ and setting in [13.1] $L_1 = l \in \mathbb{Q}_{\neq 0}$ and

[14.1]    $N = n \in \mathbb{Q}_{\neq 0}$,

with $n \neq \frac{l^2}{2}$ and $n \neq \frac{l^2}{4}$ for excluding the trivial cases, we obtain the rational solution for $M$ and $X$

[14.2]  $M = 2n - \frac{l^2}{2}$

[14.3]  $X = \pm l$

Therefore, the rational solutions for the cofactor of $X$ are

[14.4] $\frac{M}{X} = \pm \left( \frac{2n}{l} - \frac{l}{2} \right)$

and the sums of factors of $M$ (and $N$) are

[14.5]  $Z = \pm \left( \frac{2n}{l} + \frac{l}{2} \right)$

Finally, after replacing [12.2] in [5], we obtain the rational solutions for $Y$ e his cofactor $\frac{N}{Y}$

[14.6]    $Y_+ = \pm \frac{2n}{l}$,

[14.7]    $\frac{N}{Y_+} = \pm \frac{l}{2}$

and

[14.8]    $Y_- = \pm \frac{l}{2}$,

[14.9]    $\frac{N}{Y_-} = \pm \frac{2n}{l}$

We note that $Y_+ = \frac{N}{Y_-}$   $Y_- = \frac{N}{Y_+}$; therefore, the pair of solutions [14.6]-[14.7] is not essentially different from the pair [14.8]-[14.9], apart from the fact that each factor of $Y$ is switched with his cofactor[13]. Then we can arbitrarily choose the pair [14.6]-[14.7].

Therefore, the equation [1] has the following non-trivial solutions

$$M = 2n - \frac{l^2}{2}, \qquad N = n, \quad X = \pm l, \quad Y = \pm \frac{2n}{l},$$

and the consequent solutions for the cofactors of $X$ and $Y$

$$\frac{M}{X} = \pm\left(\frac{2n}{l} - \frac{l}{2}\right), \quad \frac{N}{Y} = \pm\frac{l}{2}$$

so that the sums of factors are $Z = \pm\left(\frac{2n}{l} + \frac{l}{2}\right)$ and for each arbitrary choice (except for the above restrictions) of the parameters $l$ e $n$, the pairs of rational numbers represented by the [14.1] and [14.2] has <u>two</u> IR.

---

**Case 2 ($\mathbf{a = -b}$)**

---

The [11.1] and [11.2] become respectively

[15.1]  $S = -(2N - M)$

[15.2]  $T = -M$

so that the [9.1] and [9.2] are reduced, respectively, to

[16.1]  $L_1^2 = 0 \;\Rightarrow\; X = 0$

[16.2]  $L_2^2 = 2(2N - M) \;\Rightarrow\; X = \pm L_2$

Excluding the singular case $X = 0$ and setting in [16.2] $L_2 = l \in \mathbb{Q}_{\neq 0}$ and $N = n \in \mathbb{Q}_{\neq 0}$ we obtain the same solutions as in **Case 1**.

---

**Case 3 ($\mathbf{a \neq \pm b}$)**

---

The [9.1] and [9.2] become respectively

[17.1]  $L_1^2 = \frac{a+b}{a}\left(\frac{a+b}{b}N - M\right) \;\Rightarrow\; X = \pm L_1$

[17.2]  $L_2^2 = -\frac{a-b}{a}\left(\frac{a-b}{b}N + M\right) \;\Longrightarrow\; X = \pm L_2$

In this case, the necessary condition for having rational solutions for the variable $X$ is the existence of rational solutions alternatively:

---

[13] Here if a rational number $m$ has the representation $m = ab$ then $a$ is, by definition, the cofactor of $b$ and vice-versa.

- for both the equation [17.1]-[17.2] (**Case** 3.1),

- only for the equation [17.1] (**Case** 3.2)

- only for the equation [17.2] (**Case** 3.3).

Let's examine these three sub-cases separately.

---

**Case** 3.1

---

Solve the system [17.1]-[17.2] with respect to the variables $M$ and $N$; after setting $L_1 = l_1 \in \mathbb{Q}_{\neq 0}$, $L_2 = l_2 \in \mathbb{Q}_{\neq 0}$, we obtain the solutions

[18.1] $\quad M = -\frac{1}{2}\left(\frac{a-b}{a+b}\, l_1^2 + \frac{a+b}{a-b}\, l_2^2\right)$

[18.2] $\quad N = \frac{b}{2}\left(\frac{l_1^2}{a+b} - \frac{l_2^2}{a-b}\right)$

From [7], we obtain the solutions for the variable $X$

[18.3] $\quad X = \pm l_1$

[18.4] $\quad X = \pm l_2.$

Therefore, the solutions for the cofactor of $X$ corresponding to $X = \pm l_1$ are

[18.5] $\quad \frac{M}{X} = \mp \frac{1}{2l_1}\left(\frac{a-b}{a+b}\, l_1^2 + \frac{a+b}{a-b}\, l_2^2\right)$

while those for the cofactor of $X$ corresponding to $X = \pm l_2$ are

[18.6] $\quad \frac{M}{X} = \mp \frac{1}{2l_2}\left(\frac{a-b}{a+b}\, l_1^2 + \frac{a+b}{a-b}\, l_2^2\right)$

The sums of factors of $M$ (and $N$) corresponding to $X = \pm l_1$ are

[18.7] $\quad Z = \pm \frac{(a^2 - 3b^2 + 2ab)l_1^2 - (a+b)^2 l_2^2}{2(a^2 - b^2)l_1}$

while those corresponding to $X = \pm l_2$ are

[18.8] $\quad Z = \pm \frac{(a^2 - 3b^2 - 2ab)l_2^2 - (a-b)^2 l_1^2}{2(a^2 - b^2)l_2}$

From [11.2] and [18.1] we obtain

$M + T = -\frac{a+b}{a-b}\, l_2^2$

$M - T = -\frac{a-b}{a+b}\, l_1^2$

By replacing these last two expressions in [5] we obtain the rational solutions for $Y$ and for the cofactor $\frac{N}{Y}$ corresponding to $X = \pm l_1$

[18.9] $\qquad Y_+ = \pm \frac{(a-b)l_1^2-(a+b)l_2^2}{2(a-b)l_1}$

[18.10] $\qquad \frac{N}{Y_+} = \pm \frac{bl_1}{a+b}$

[18.11] $\qquad Y_- = \pm \frac{bl_1}{a+b}$

[18.12] $\qquad \frac{N}{Y_-} = \pm \frac{(a-b)l_1^2-(a+b)l_2^2}{2(a-b)l_1}$

and to $X = \pm l_2$

[18.13] $\qquad Y_+ = \mp \frac{bl_2}{a-b}$

[18.14] $\qquad \frac{N}{Y_+} = \pm \frac{(a+b)l_2^2-(a-b)l_1^2}{2(a+b)l_2}$

[18.15] $\qquad Y_- = \pm \frac{(a+b)l_2^2-(a-b)l_1^2}{2(a+b)l_2}$

[18.16] $\qquad \frac{N}{Y_-} = \mp \frac{bl_2}{a-b}$

Again we note that $Y_+ = \frac{N}{Y_-}$ e $Y_- = \frac{N}{Y_+}$; therefore, the pair of solutions [18.9]-[18.10] corresponding to $X = \pm l_1$ is not essentially different from the pair [18.11]-[18.12], apart from the fact that each factor of $Y$ is switched with his cofactor. Then we can arbitrarily choose the pair [18.9]-[18.10]. The same property obviously holds for the solutions corresponding to $X = \pm l_2$, of which we choose the pair [18.13]-[18.14].

For excluding the trivial solutions, we set in addition $\left(\frac{l_1}{l_2}\right)^2 \neq \pm \frac{a+b}{a-b}$.

Now we have the following sub-cases that we examine separately:

- $l_1 \neq l_2 \in \mathbb{Q}_{\neq 0}$ $\qquad$ (*Case* 3.1.1)

- $l_1 = l_2 = l \in \mathbb{Q}_{\neq 0}$ $\quad$ (*Case* 3.1.2)

---

| ***Case* 3.1.1 ($a \neq \pm b$, $\qquad l_1 \neq l_2$)** |
|---|

In this sub-case for each arbitrary choice (except for the above restrictions) of the parameters $a$, $b$, $l_1$ e $l_2$, the pairs of rational numbers represented by the [18.1] and [18.2] have <u>four</u> IR; the respective factors, in turn, can be represented by the formulas from [18.3] to [18.6] (for the factors of $M$) and the formulas [18.9]-[18.10] and [18.13]-[18.14] (for the factors of $N$).

---

**Case** 3.1.2 $((a \neq \pm b, l_1 = l_2 = l)$

---

In this sub-case the [18.1] and [18.2] become respectively

[19.1] $\quad M = -\frac{a^2+b^2}{a^2-b^2}l^2$

[19.2] $\quad N = -\frac{b^2}{a^2-b^2}l^2$

Thus, each pair of rational numbers representable by [19.1] and [19.2] has <u>two</u> IR, whose factors are

[19.3] $\quad X = \pm l$

[19.4] $\quad \frac{M}{X} = \mp\frac{a^2+b^2}{a^2-b^2}l$

[19.5] $\quad Y = \mp\frac{b}{a-b}l$

[19.6] $\quad \frac{N}{Y} = \pm\frac{b}{a+b}l$

and the sums of factors of $M$ (and $N$) are

[19.7] $\quad Z = \mp\frac{2b^2}{a^2-b^2}l$

---

**Case** 3.2

---

In [17.1] we set $L_1 = l \in \mathbb{Q}_{\neq 0}$ and

[20.1] $\quad N = n \in \mathbb{Q}_{\neq 0},$

with $n \neq \frac{abl}{(a+b)^2}$ for excluding the trivial solutions, and we obtain the rational solutions for $M$ and $X$

[20.2] $\quad M = \frac{a+b}{b}n - \frac{al^2}{a+b}$

[20.3] $\quad X = \pm l$

Therefore, the rational solutions for the cofactor of $X$ are

[20.4] $\quad \frac{M}{X} = \pm\left(\frac{a+b}{bl}n - \frac{al}{a+b}\right)$

and the sums of factors of $M$ (and $N$) are

[20.5] $\quad Z = \pm\frac{b^2l^2+(a+b)^2n}{bl(a+b)}$

From [11.2] and [20.2] we obtain

$$M + T = \frac{2(a+b)}{b}n - l^2$$

$$M - T = -\frac{a-b}{a+b}l^2$$

By replacing these last two expressions in [5] we obtain the pairs of rational solutions for $Y$ and the cofactor $\frac{N}{Y}$

[20.6]    $Y_+ = \pm\frac{a+b}{bl}n,$

[20.7]    $\frac{N}{Y_+} = \pm\frac{bl}{a+b}$

and

[20.8]    $Y_- = \pm\frac{bl}{a+b},$

[20.9]    $\frac{N}{Y_-} = \pm\frac{a+b}{bl}n$

Like before $Y_+ = \frac{N}{Y_-}$ and $Y_- = \frac{N}{Y_+}$; therefore the pair of solutions [20.6]-[20.7] is not essentially different from the pair [20.8]-[20.9], apart from the fact that each factor of $Y$ is switched with the respective cofactor. Then we can arbitrarily choose the pair [20.6]-[20.7].

Thus, the equation [1] has the following non-trivial solutions

$$M = \frac{a+b}{b}n - \frac{al^2}{a+b}, \ \ N = n, \ \ \ X = \pm l, \ \ \ Y = \pm\frac{a+b}{bl}n$$

and the consequent solutions for the cofactors of $X$ e $Y$

$$\frac{M}{X} = \frac{a+b}{bl}n - \frac{al}{a+b}, \ \ \ \frac{N}{Y} = \pm\frac{bl}{a+b},$$

so that, the sum of factors are $Z = \pm\frac{b^2 l^2 + (a+b)^2 n}{bl(a+b)}$ and for each arbitrary choice (except for the above restrictions) of the parameters $a, b, l$ e $n$, the pairs of rational numbers represented by the [20.1] and [20.2] has <u>two</u> IR[14].

---

**Case** 3.3

In [17.2] we set $L_2 = l \in \mathbb{Q}_{\neq 0}$ and

[21.1]    $N = n,$

---

[14] If we set $a = b$ these formulas coincide with those of **Case** 1.
Francesco Trimarchi, Rational points on elliptic curves and representations of rational numbers as the product of two rational factors, Milano, December 2018 with $n \neq -\frac{abl^2}{(a-b)^2}$ for excluding trivial solutions, and we obtain the rational solutions for $M$ and $X$

[21.2] $\qquad M = -\left(\frac{a-b}{b}n + \frac{al^2}{a-b}\right)$

[21.3] $\qquad X = \pm l$

Therefore, the rational solutions for the cofactor of $X$ are

[21.4] $\qquad \frac{M}{X} = \mp \left(\frac{a-b}{bl}n + \frac{al}{a-b}\right)$

and the sums of factors of $M$ (and $N$) are

[21.5] $\quad Z = \mp \frac{b^2l^2+(a-b)^2n}{bl(a-b)}$

From [11.2] and [21.2] we obtain

$$M + T = -\frac{a+b}{a-b}l^2$$

$$M - T = -\left(\frac{2(a-b)}{b}n + l^2\right)$$

By replacing these last two expressions in [5] we obtain the pairs of rational solutions for $Y$ and the cofactor $\frac{N}{Y}$

[21.6] $\qquad Y_+ = \mp \frac{bl}{a-b}$,

[21.7] $\qquad \frac{N}{Y_+} = \mp \frac{a-b}{bl}n$

and

[21.8] $\qquad Y_- = \mp \frac{a-b}{bl}n$,

[21.9] $\qquad \frac{N}{Y_-} = \mp \frac{bl}{a-b}$

Like before, $Y_+ = \frac{N}{Y_-}$ and $Y_- = \frac{N}{Y_+}$; then the pair of solutions [21.6]-[21.7] is not essentially different from the pair [21.8]-[21.9], apart from the fact that each factor of $Y$ is switched with the respective cofactor. Then we can arbitrarily choose the pair [21.6]-[21.7].

Thus, the equation [2] has the following solutions

$$M = -\left(\frac{a-b}{b}n + \frac{al^2}{a-b}\right), \qquad N = n, \qquad X = \pm l, \qquad Y = \mp \frac{bl}{a-b}$$

and the consequent solutions for the cofactor of $X$ and $Y$

$$\frac{M}{X} = \mp \left( \frac{a-b}{bl} n + \frac{al}{a-b} \right), \quad \frac{N}{Y} = \mp \frac{a-b}{bl} n$$

so that the sums of factors are $Z = \mp \frac{b^2 l^2 + (a-b)^2 n}{bl(a-b)}$ and, for each arbitrary choice (except for the above restrictions) of the parameters $a, b, l$ e $n$, the pairs of rational numbers represented by the [21.1] and [21.2] have <u>two</u> iso-additive representations[15].

<p style="text-align:center">*  *  *</p>

Since in all possible cases $\rho(m, n) \leq 4$ the IRT is proved[16] ∎

## §3 *Rational points on elliptic curves and IRT*

We consider the elliptic curve

[1]    $E: W^2 = U(U-m)(U-n)$,

where $m \neq n \in \mathbb{Q}_{\neq 0}$, and let $E(\mathbb{Q})$ the group of rational points on $E$. We prove the following:

**Theorem 1** (*Iso − additive Rational Points Theorem − IRPT*) *If there exist four non-zero rational numbers* $m_1$, $m_2$, $n_1$ *and* $n_2$ *such that* $m_1 m_2 = m$, $n_1 n_2 = n$, $m_1 m_2 \neq n_1 n_2$ *and* $m_1 + m_2 = n_1 + n_2$, *then the set* $\mathcal{P}(m,n) = \{\pm P_j, \ j = 1,2,3,4\}$ *where*

$P_1 = (m_1 n_1, \ (m_2 - n_1) m_1 n_1), \quad -P_1 = (m_1 n_1, \ -(m_2 - n_1) m_1 n_1)$

$P_2 = (m_1 n_2, \ (m_1 - n_1) m_1 n_2), \quad -P_3 = (m_1 n_2, \ -(m_1 - n_1) m_1 n_2)$

$P_3 = (m_2 n_1, \ (m_1 - n_1) m_2 n_1), \quad -P_3 = (m_2 n_1, \ -(m_1 - n_1) m_2 n_1)$

$P_4 = (m_2 n_2, \ (m_2 - n_1) m_2 n_2), \quad -P_4 = (m_2 n_2, \ -(m_2 - n_1) m_2 n_1)$

is a sub-set of $E(\mathbb{Q})$.

Proof. The statement follows immediately from the definition of these points (replace their coordinates in the equation [1] and check the result).

Anyway, we give an equivalent proof that relates the set $\wp(m,n) = \{\pm p_j, \ j = 1,2,3,4\}$ as defined in Ch. 2 to $E(\mathbb{Q})$.

---

[15] Also in this case if we set $a = b$ these formulas coincide with those of **Case** 1.

[16] As we showed in note 10, when $\rho(m,n) = 2$ we can rewrite all these solutions in the form $m = t^2 - r^2$ and $n = t^2 - s^2$ so that $m_1 = \pm(t+r)$, $m_2 = \pm(t-r)$, $n_1 = \pm(t+s)$ and $n_2 = \pm(t-s)$. When $\rho(m,n) = 4$, in addition we have $m = t'^2 - r'^2$ and $n = t'^2 - s'^2$ so that $m_1 = \pm(t'+r')$, $m_2 = \pm(t'-r')$, $n_1 = \pm(t'+s')$ and $n_2 = \pm(t'-s')$.

After some manipulation[17] the equation $X + \frac{m}{X} = Y + \frac{m}{Y}$ become

[2]     $W^2 = U(U - m)(U - n)$

where $U = XY$ and $W = \left(\frac{n}{Y} - X\right)XY$, or equivalently $\left(\frac{m}{X} - Y\right)XY \to W$.

Therefore, by the transformation $XY \to U$ and $\left(\frac{m}{X} - Y\right)XY \to W$, to the solution $p_1 = (x, \ y)$ of the equation $X + \frac{m}{X} = Y + \frac{m}{Y}$, where $x \neq y \in \mathbb{Q}_{\neq 0}$, corresponds on the elliptic curve $W^2 = U(U - m)(U - n)$ the rational point $P_1 = (u_1, w_1)$, where $u_1 = xy$ and $w_1 = \left(\frac{m}{x} - y\right)xy$.

On the other hand, in Ch. 2 we have seen that not only $p_1$ but any element of the set $\wp(m,n) = \{\pm p_j, \ j = 1,2,3,4\}$, generated from $p_1$ by switching factors and changing sign, satisfies the equation $X + \frac{m}{X} = Y + \frac{n}{Y}$ so that to any solution $\pm p_j \in \wp(m,n)$ of this equation corresponds on $E$ a rational point, whose coordinates are given by the transformation $XY \to U$ and $\left(\frac{m}{X} - Y\right)XY \to W$, or equivalently $\left(\frac{n}{Y} - X\right)XY \to W$.

Therefore, this transformation defines a map from the set $\wp(m,n)$ and a set, say $\mathcal{P}(m,n) \subseteq E(\mathbb{Q})$ whose elements are

$P_1 = \left(xy, \ \ \left(\frac{m}{x} - y\right)xy\right)$

$P_2 = \left(x\frac{n}{y}, \ \ (x - y)x\frac{n}{y}\right)$

$P_3 = \left(\frac{m}{x}y, \ \ (x - y)\frac{m}{x}y\right)$

$P_4 = \left(\frac{m}{x}\frac{n}{y}, \ \ \left(\frac{m}{x} - y\right)\frac{m}{x}\frac{n}{y}\right)$

and their opposite $-P_j, \ j = 1,2,3,4$.

For homogenize the symbolism let $x = m_1$, $\frac{m}{x} = m_2$, $y = n_1$ and $\frac{n}{y} = n_2$ so that we obtain the same formulas of the statement∎

As we have seen in Ch. 2, for given $m$ and $n$ when $\rho(m,n) = 2$ the only rational solutions of the equation $X + \frac{m}{X} = Y + \frac{n}{Y}$ are the elements of $\wp(m,n)$ so that $\Gamma(m,n) = \wp(m,n)$; when $\rho(m,n) = 4$ $\Gamma(m,n) = \wp(m,n) \cup q(m,n)$, where the set $q_i(m,n) = \{\pm q_j, \; j = 1,2,3,4\}$ has already been defined in Ch. 2; this define a second map $q_i(m,n) \to \mathcal{Q}(m,n) \subseteq E(\mathbb{Q})$ and the points of $\mathcal{Q}(m,n)$ are

$$Q_1 = \left(x'y', \quad \left(\frac{m}{x'} - y'\right)x'y'\right)$$

$$Q_2 = \left(x'\frac{n}{y'}, \quad (x' - y')x'\frac{n}{y'}\right)$$

$$Q_3 = \left(\frac{m}{x'}y', \; (x' - y')\frac{m}{x'}y'\right)$$

$$Q_4 = \left(\frac{m}{x'}\frac{n}{y'}, \; \left(\frac{m}{x'} - y'\right)\frac{m}{x'}\frac{n}{y'}\right)$$

and their opposite $-Q_j$, $j = 1,2,3,4$.

---

[17] In Introduction we showed that the equation $X + \frac{m}{X} = Y + \frac{n}{Y}$ is equivalent to $\frac{X}{Y} = \frac{XY - m}{XY - n}$; (see note ); then, if we multiply each side for $XY(XY - n)^2$ we easily obtain the equation $\left[\left(\frac{m}{X} - Y\right)XY\right]^2 = XY(XY - m)(XY - n)$; finally, by setting $U = XY$ and $W = \left(\frac{m}{X} - Y\right)XY$, or equivalently $W = \left(\frac{n}{Y} - X\right)XY$, we obtain the equation $W^2 = U(U - m)(U - n)$.

Now we show the main properties of elliptic curves of this type when $\rho(m,n) = 2$.

Let $\mathcal{T} = (0,0)$, $\mathcal{M} = (m,0)$ and $\mathcal{N} = (n,0)$ the points of order two and $\mathcal{O}$ the point at infinity. We use the symbol "+" for the addition of points on $E(\mathbb{Q})$ and the expression $[k]P$ for the addition $k$-times of the same point $P \in E(\mathbb{Q})$; moreover, given three points on $E(\mathbb{Q})$ such that $P + Q = R$ (i.e. $P$, $Q$ and $-R$ are collinear), let $\lambda$ and $\nu$ respectively the slope and the intercept of the straight line $w = \nu + \lambda u$ which intersects $E$ in $P$, $Q$ and $-R$. Moreover, points of order two satisfy the obvious identities $\mathcal{T} = -\mathcal{T}$, $\mathcal{M} = -\mathcal{M}$ and $\mathcal{N} = -\mathcal{N}$, as well the identity $\mathcal{T} + \mathcal{M} + \mathcal{N} = \mathcal{O}$ (the sum of any two points of order two gives the third, because they are collinear, and the straight line passing through them coincides with the $U$-axis).

Preliminary, we prove two Lemmas concerning the sub-set $\mathcal{P}(m,n) \cup \{\mathcal{T}, \mathcal{M}, \mathcal{N}, \mathcal{O}\}$ which are useful for simplifying the proof of some proposition about the points of finite order on the curve $W^2 = U(U - m)(U - n)$, where $m \neq n \in \mathbb{Q}_{\neq 0}$ and $(m,n)$ have two IR.

Obviously, the same properties stated in the following Lemmas hold when roots of the cubic have four IR and $E(\mathbb{Q})$ contains also the sub-set $\mathcal{Q}(m,n)$.

*Lemma* 1 (*addition*). *The points of* $\mathcal{P}(m,n) \cup \{\mathcal{T}, \mathcal{M}, \mathcal{N}, \mathcal{O}\}$ *satisfy the following identities (in bracket the coefficients of the straight line passing through the three points involved):*

$$P_1 + P_4 = \mathcal{T}, \qquad \left(\lambda = \frac{m}{x} - y, \quad \nu = 0\right),$$

$$P_2 + P_3 = \mathcal{T}, \qquad (\lambda = x - y, \quad \nu = 0),$$

$$P_1 + P_3 = \mathcal{M}, \qquad \left(\lambda = -y, \qquad \nu = y\frac{m}{x}\right),$$

$$P_2 + P_4 = \mathcal{M}, \qquad \left(\lambda = \frac{m}{x} - y, \quad \nu = -\frac{mn}{y}\right),$$

$$P_1 - P_2 = \mathcal{N}, \qquad \left(\lambda = -x, \qquad \nu = x\frac{n}{y}\right),$$

$$P_3 - P_4 = \mathcal{N}, \qquad \left(\lambda = -\frac{m}{x}, \qquad \nu = -\frac{mn}{x}\right)$$

```
Proof. We get these identities by applying the well-known formulas for the
sum of two rational points on elliptic curves[18]∎
```

*Lemma* 2 (*duplication*)

$$[2]P_j = \left(\left(\frac{z}{2}\right)^2, \ t_j \frac{z}{2}\left(\frac{z}{2} - x\right)\left(\frac{z}{2} - y\right)\right), \qquad t_j = \begin{cases} -1, & j = 1,2 \\ 1, & j = 3,4 \end{cases}$$

where $z = x + \frac{m}{x}\left(= y + \frac{n}{y}\right)$ and $x \neq y \in \mathbb{Q}_{\neq 0}$.

```
Proof. Let U(Pj) = U(−Pj), j = 1,2,3,4 the U−coordinates of each point of P(ℚ);
preliminary, note that ∑⁴ⱼ₌₁ U(Pj) = z² and consequently
```

$$\left(\frac{z}{2}\right)^2 \left[\left(\frac{z}{2}\right)^2 - m\right]\left[\left(\frac{z}{2}\right)^2 - n\right] = \left[\frac{z}{2}\left(\frac{z}{2} - x\right)\left(\frac{z}{2} - y\right)\right]^2$$

```
Therefore, the points P = {(z/2)², z/2(z/2 − x)(z/2 − y)} and −P belong to E(ℚ).
```

```
On the other hand, applying duplication formula[19] for any point Pj =
```
$(u_j, \ w_j) \in \mathcal{P}(\mathbb{Q}),$ we have $U([2]P_j) = \left(\frac{u_j^2 - mn}{2w_j}\right)^2 = U(P)$ and consequently the $W-$coordinates of $[2]P_j$ are $W([2]P_j) = t_j W(P)$, where $t = \begin{cases} -1, & j = 1,2 \\ 1, & j = 3,4 \end{cases}$

```
This complete the proof ∎
```

```
From Lemma 1 we obtain the three following corollaries:
```

*Corollary* 1. *Given a fixed point of* $\mathcal{P}(m,n)$, *any point of* $\mathcal{P}(m,n)$ *can be represented as the sum of this fixed point, or its opposite, and a suitable point of order two, or the point at infinity, in this way:*

*Fixed point*                 *Representations of* $P_j$

| | | | | |
|---|---|---|---|---|
| $P_1$ | $P_1 = P_1 + \mathcal{O},$ | $P_2 = P_1 + \mathcal{N},$ | $P_3 = -P_1 + \mathcal{M},$ | $P_4 = -P_1 + \mathcal{T}$ |
| $P_2$ | $P_1 = P_2 + \mathcal{N},$ | $P_2 = P_2 + \mathcal{O},$ | $P_3 = -P_2 + \mathcal{T},$ | $P_4 = -P_2 + \mathcal{M}$ |
| $P_3$ | $P_1 = -P_3 + \mathcal{M},$ | $P_2 = -P_3 + \mathcal{T},$ | $P_3 = P_3 + \mathcal{O},$ | $P_4 = P_3 + \mathcal{N}$ |

$$P_4 \qquad P_1 = -P_4 + \mathcal{T}, \qquad P_2 = -P_4 + \mathcal{M}, \quad P_3 = P_4 + \mathcal{N}, \qquad P_4 = P_4 + \mathcal{O}$$

Proof. These identities follow immediately from the six identities in *Lemma* 1 and the properties of the point at infinity (for simplicity we omit the obvious representations of $-P_j, j = 1,2,3,4$) ∎

*Corollary* 2. *Let* $P_1$ *a fixed point. The following identities holds*[20]

$$[2k]P_1 = [2k]P_2 = -[2k]P_3 = -[2k]P_4$$

$$[2k + 1]P_1 = [2k + 1]P_1 + \mathcal{O},$$

$$[2k + 1]P_2 = [2k + 1]P_1 + \mathcal{N}$$

$$[2k + 1]P_3 = -[2k + 1]P_1 + \mathcal{M}$$

$$[2k + 1]P_4 = -[2k + 1]P_1 + \mathcal{T}$$

Proof. The statement follows immediately adding $2k$-times (or $2k + 1$-times) each member of the identities showed in first row of *Corollary* 1 and from the fact that the addition $h$-times of a point of order two gives the point itself, if $h = 2k + 1$, or the point at infinity, if $h = 2k$ ∎

*Corollary* 3. *The following identities holds*

$$P_1 + P_2 = \mathcal{N} + [2]P_1$$

$$P_3 + P_4 = \mathcal{N} - [2]P_1$$

$$P_1 - P_3 = \mathcal{M} + [2]P_1$$

$$P_2 - P_4 = \mathcal{M} + [2]P_1$$

$$P_1 - P_4 = \mathcal{T} + [2]P_1$$

$$P_2 - P_3 = \mathcal{T} + [2]P_1$$

Proof. We obtain these identities from *Corollary* 1 and the property $[2]P_1 = [2]P_2 = -[2]P_3 = -[2]P_4$ which is a special case of the property $[2k]P_1 = [2k]P_2 = -[2k]P_3 = -[2k]P_4$ proved in *Corollary* 2 ∎

The six identities in *Lemma* 1, together the six in *Corollary* 3 and the properties of $\{\mathcal{T}, \mathcal{M}, \mathcal{N}, \mathcal{O}\}$, give the addition of any pair of points of the set $\mathcal{P}(m,n) \cup \{\mathcal{T}, \mathcal{M}, \mathcal{N}, \mathcal{O}\}$.

While any addition defined by first six identities gives a third point which still belongs to $\mathcal{P}(m,n) \cup \{\mathcal{T}, \mathcal{M}, \mathcal{N}, \mathcal{O}\}$, this could not be the case for the additions defined by the remaining six identities, where the third point is the sum of a fixed point ($[2]P_1$ or its opposite) and a suitable point of order two.

---

[20] It's easy to obtain similar identities if the fixed point is $P_2$, $P_3$ or $P_4$, as in *Corollary* 1.

However, under some conditions also the additions defined by the remaining six identities give a point of $\mathcal{P}(m\,n) \cup \{\mathcal{T}, \mathcal{M}, \mathcal{N}, \mathcal{O}\}$; when it happens, this set is close under addition and consequently it is a group.

We find some of these conditions showing what happens when each point

$$\mathcal{N} + 2P_1, \quad \mathcal{N} - 2P_1, \quad \mathcal{M} + 2P_1, \quad \mathcal{T} + 2P_1$$

coincide with a point of the set $\{\mathcal{T}, \mathcal{M}, \mathcal{N}, \mathcal{O}\}$ which ensures the closure under addition of the set $\mathcal{P}(m,n) \cup \{\mathcal{T}, \mathcal{M}, \mathcal{N}, \mathcal{O}\}$.

This is a very interesting case, although a complete assessment of these conditions should require verifying what happens also when each of these points coincide whit a point of $\mathcal{P}(m,n)$. Here we prove the following

*Theorem* 2. *If one of the following identities holds*

$$\mathcal{N} + 2P_1 = \mathcal{M}, \ \mathcal{N} + 2P_1 = \mathcal{T}, \ \mathcal{N} + 2P_1 = \mathcal{O}$$

*then*

  i)    $\mathcal{P}(m,n) \cup \{\mathcal{T}, \mathcal{M}, \mathcal{N}, \mathcal{O}\}$ *is a group;*

  ii)   $\#\mathcal{P}(m,n) = 4$, *i.e. the points of $\mathcal{P}(m,n)$ coincide two-by-two;*

  iii)  *any point of $\mathcal{P}(m,n)$ has order four.*

Proof. Preliminary, note that if $\mathcal{N} + 2P_1 = \mathcal{N}$ (so that also $\mathcal{N} - 2P_1 = \mathcal{N}$) then $[2]P_1 = \mathcal{O}$, i.e $P_1$ would be a point of order two, which contradicts our assumptions[21].

Then we start by setting $\mathcal{N} + 2P_1 = \mathcal{M}$ so that also $\mathcal{N} - 2P_1 = \mathcal{M}$; this implies $2P_1 = \mathcal{T}$ and consequently $\mathcal{M} + 2P_1 = \mathcal{N}$ and $\mathcal{T} + 2P_1 = \mathcal{O}$. Therefore, the identities of *Corollary* 3 become

$$P_1 + P_2 = \mathcal{M}, \qquad P_1 - P_3 = \mathcal{N}, \qquad P_1 - P_4 = \mathcal{O},$$

$$P_3 + P_4 = \mathcal{M}, \qquad P_2 - P_4 = \mathcal{N}, \qquad P_2 - P_3 = \mathcal{O}.$$

and ensure the closure under addition of $\mathcal{P}(m,n) \cup \{\mathcal{T}, \mathcal{M}, \mathcal{N}, \mathcal{O}\}$. Moreover, $P_4 = P_1$ and $P_3 = P_2$ so that $\mathcal{P}(m,n) = \{P_1, -P_1, P_2, -P_2\}$ and $\#\mathcal{P}(\mathbb{Q}) = 4$.

Finally, from the identity $[2]P_1 = \mathcal{T}$ we get $[3]P_1 = \mathcal{T} + P_1 = -P_1$ and $[4]P_1 = \mathcal{O}$. Since $P_1$, $[2]P_1$ and $[3]P_1$ do not coincide with the point at infinity, $P_1$

---

[21] If $[2]P_1$ should coincide whith $\mathcal{O}$ then $P_1$ should be equal to its opposite. Provided that by assumption $x$ and $y$ must non-zero and distinct, $P_1 = -P_1$ holds if and only if $\frac{m}{x} - y = \frac{n}{y} - x = 0$ or equivalently $m = n = xy$, which contradicts the assumption that $m$ and $n$ are distinct. Consequently $P_1$ cannot coincide with a point of order two (or less since obviously $P_1$ can not be coincide with the point at infinity)

as well $-P_1$ are of order four; on the other hand, from *Corollary* 2 $[2]P_1 = [2]P_2$ therefore $[2]P_2 = \mathcal{T}$, $[3]P_2 = \mathcal{T} + P_2 = -P_2$ and $[4]P_2 = \mathcal{O}$ so that also $P_2$ as well $-P_2$ are of order four.

Now we repeat the previous step setting either $\mathcal{N} + [2]P_1 = \mathcal{T}$ or $\mathcal{N} + [2]P_1 = \mathcal{O}$; at the end we obtain

| | $\mathcal{P}(m,n)$ | $[2]P_1$ | *Additions* |
|---|---|---|---|
| $\mathcal{N} + 2P_1 = \mathcal{M}$ | $\{P_1, -P_1, P_2, -P_2\}$ | $\mathcal{T}$ | $P_1 + P_2 = \mathcal{M}$, $P_1 - P_2 = \mathcal{N}$ |
| $\mathcal{N} + 2P_1 = \mathcal{T}$ | $\{P_1, -P_1, P_2, -P_2\}$ | $\mathcal{M}$ | $P_1 + P_2 = \mathcal{T}$, $P_1 - P_2 = \mathcal{N}$ |
| $\mathcal{N} + 2P_1 = \mathcal{O}$ | $\{P_1, -P_1, P_3, -P_3\}$ | $\mathcal{N}$ | $P_1 + P_3 = \mathcal{M}$, $P_1 - P_3 = \mathcal{T}$ |

This complete the proof ∎

Now, from *Lemma* 2 (*duplication*) we obtain the explicit formulas for the parameters of the curve $W^2 = U(U - m)(U - n)$ and the coordinates of each points of $\mathcal{P}(m,n) \cup \{\mathcal{T}, \mathcal{M}, \mathcal{N}, \mathcal{O}\}$, setting $[2]P_1$ alternatively equal to $\mathcal{T} = (0,0)$, $\mathcal{M} = (m,0)$ or $\mathcal{N} = (n,0)$ and remembering that: i) $x$ and $y$ must be non-zero and distinct; ii) $z = x + \frac{m}{x}\left(= y + \frac{n}{y}\right)$.

---

**Case 1**: $[2]P_1 = \mathcal{T}$

The equation $\left\{\left(\frac{z}{2}\right)^2, \ \frac{z}{2}\left(\frac{z}{2} - x\right)\left(\frac{z}{2} - y\right)\right\} = (0,0)$, is satisfied when $z = 0$, $x \neq \frac{z}{2}$ and $y \neq \frac{z}{2}$ (both the last two inequalities are necessary because $x$ and $y$ cannot be zero by assumption). Therefore, equating to zero the sum of factors of $m$ and $n$ we obtain the IR $m = -x^2 = x \cdot (-x)$ and $n = -y^2 = y \cdot (-y)$; consequently $\mathcal{M} = (-x^2, 0)$ and $\mathcal{N} = (-y^2, 0)$. This Case corresponds to the elliptic curve $W^2 = U(U + x^2)(U + y^2)$ such that the non-zero roots of the cubic are the opposite of a square.

From *Theorem* 2 the points of $\mathcal{P}(m,n)$ are

$P_1 = \{xy, \ -(x + y)xy\}$

$P_2 = \{-xy, \ (x + y)xy\}$

and their opposite.

---

**Case 2**: $[2]P_1 = \mathcal{M}$

The equation $\left\{\left(\frac{z}{2}\right)^2, \ \frac{z}{2}\left(\frac{z}{2} - x\right)\left(\frac{z}{2} - y\right)\right\} = (m, 0)$ is satisfied when $\left(\frac{z}{2}\right)^2 = m$ and $x = \frac{z}{2}$ but $y \neq \frac{z}{2}$ or vice-versa (by assumption, $x$ and $y$ cannot be equal). Now $x = \frac{z}{2}$ implies $z = 2x$; by substituting $2x$ for $z$ in the equation $z = x + \frac{m}{x}\left(= y + \frac{n}{y}\right)$ we obtain $m = x^2$ and $n = x^2 - (x - y)^2$; therefore, $\mathcal{M} = (x^2, 0)$, $\mathcal{N} = (x^2 - (x - y)^2, 0)$ and $(m, n)$ have the IR $m = x \cdot x$ and $n = y \cdot (2x - y)$ such that the sum of factors is $2x$. This Case corresponds to the elliptic curve $W^2 = U(U - x^2)\left(U - (x^2 - (x - y)^2)\right)$ such that one of the non-zero roots of the cubic is a square while the other is the difference between the first root and some square.

From *Theorem* 1 the points of $\mathcal{P}(m, n)$ are now

$P_1 = \{xy, \ (x - y)xy\}$

$P_2 = \{x(2x - y), \quad (x - y)x(2x - y)\}$

and their opposite.

If we set $\left(\frac{z}{2}\right)^2 = m$, as before, but $y = \frac{z}{2}$ and $x \neq \frac{z}{2}$ (instead of $x = \frac{z}{2}$ and $y \neq \frac{z}{2}$), we obtain the same results apart from a slight change of symbolism: $m = y^2 - (x - y)^2$ (instead of $x^2$) and $n = y^2$ (instead of $x^2 - (x - y)^2$). The other changes of symbolism follow immediately.

---

**Case 3:** $[2]P_1 = \mathcal{N}$

---

The equation $\left\{\left(\frac{z}{2}\right)^2, \ \frac{z}{2}\left(\frac{z}{2} - x\right)\left(\frac{z}{2} - y\right)\right\} = (n, 0)$ is satisfied when $\left(\frac{z}{2}\right)^2 = n$ and $x = \frac{z}{2}$ but $y \neq \frac{z}{2}$, or vice-versa, and this leads to the same results showed for the Case 2, apart from the consequent changes of symbolism.

These results are the basis for a complete assessment of the properties of the elliptic curves of the type

$$W^2 = U(U - x^2)\left(U - (x^2 - (x - y)^2)\right)$$

Here two examples, already introduced in Ch. 1:

1) Equation $W^2 = U(U - 1)(U + 3)$. We derive it from the equation above by setting $x = 1$ and $y = 3$ so that the two non-zero roots of the cubic become $m = 1^2 = 1$ and $n = 1^2 - (1 - 3)^2 = 1^2 - 2^2 = -3$ have the two IR $m = 1 \cdot 1 = -1 \cdot (-1)$ and $n = -1 \cdot 3 = 1 \cdot (-3)$; the points of order two are $\mathcal{T} = (0, 0)$, $\mathcal{M} = (1, 0)$ and $\mathcal{N} = (-3, 0)$ while the point of $\mathcal{P}(m, n)$ are

$P_1 = \{xy, \ (x - y)xy\} = (3, -6)$

$P_2 = \{x(2x - y), \ (x - y)x(2x - y)\} = (-1, 2)$

and their opposite, $-P_1 = (3, 6)$ and $-P_2 = (-1, -2)$; moreover, $[2]P_1 = [2]P_2 = \mathcal{M} = (1, 0)$ so that $[4]P_1 = [4]P_2 = \mathcal{O}$; therefore all points of $\mathcal{P}(m,n)$ are of order four and $\mathcal{P}(m,n) \cup \{\mathcal{T},\ \mathcal{M},\ \mathcal{N},\ \mathcal{O}\} = E(\mathbb{Q})$. Note that $\#E(\mathbb{Q}) = 8$ while the existence of four distinct square in arithmetic sequence would require an additional point on $E(\mathbb{Q})$; this solves (negatively) the question of four square in arithmetic sequence. Moreover, the pair of non-zero roots $(1,-3)$ is a solution, for the variables $M$ and $N$, of the equation $X + \frac{M}{X} = Y + \frac{N}{Y}$ (Ch. 2, p. 16) whose essentially different solutions, in parametric form, are

$$M = -\left(\frac{a-b}{b}n + \frac{al^2}{a-b}\right), \qquad N = n, \qquad X = \pm l, \qquad Y = \mp\frac{bl}{a-b}$$

Therefore, the solutions for the cofactor of $X$ and $Y$ are

$$\frac{M}{X} = \mp\left(\frac{a-b}{bl}n + \frac{al}{a-b}\right), \qquad \frac{N}{Y} = \mp\frac{a-b}{bl}n$$

while the sums of factors are $Z = \mp\frac{b^2l^2+(a-b)^2n}{bl(a-b)}$.

These solutions depend on four parameters $(a,b,l,n)$ and if we set $a = 2$, $b = 3$, $l = -1$ and $n = -3$ we exactly obtain $m = 1$ as well the IR of $m$ and $n$ showed above.

2) the special case of the Frey curve $W^2 = U(U - a^{2k})(U + b^{2k})$ under the assumption $a^{2k} + b^{2k} = c^{2k}$ for some integer $a,b,c,k > 0$. If it were true the two non-zero roots of the cubic, $a^{2k}$ and (as a consequence of our assumption) $-b^{2k} = a^{2k} - c^{2k}$, would have the IR $a^{2k} = a^k \cdot a^k$ and $-b^{2k} = (a^k + c^k) \cdot (a^k - c^k)$ such that the sum of factors is $2a^k$. Therefore, also this curve would have four points of order four which, together $\mathcal{T} = (0,0)$, $\mathcal{M} = (a^{2k},0)$, $\mathcal{N} = (a^{2k} - c^{2k},0)$ and the point at infinity, form the group of rational points on the curve. It is well-known that this curve exists only for $k = 1$ (Pytagorean case, see also Ch. 1, p. 2) because of Last Fermat Theorem[22]. Therefore, if we could prove – without using the theory of modular forms – that there is no elliptic curve of the type $W^2 = U(U - a^{2k})(U + b^{2k})$ with the properties of $E(\mathbb{Q})$ previously showed (related to IRT/IRPT), then the Last Fermat Theorem would be automatically proved for the well-known case of even exponent. However, the approach based on IRT/IRPT could be the basis for an elementary proof of the relevant case (odd prime exponent).

---

[22] Silverman-Tate, pp. 245-255.